\documentclass [reqno]{amsart}
\usepackage[english]{babel}
\usepackage{textcomp}
\usepackage[leqno]{amsmath}
\usepackage{amsfonts,amssymb,amsthm,mathrsfs,float}
\usepackage{graphicx}
\usepackage{lmodern}
\usepackage{hyperref}
\usepackage{etoolbox}
\usepackage{mathptmx}
\usepackage{mathtools} 
\usepackage{cancel}
\patchcmd{\abstract}{\scshape\abstractname}{\textbf{\abstractname}}{}{}
\hypersetup{urlcolor=blue, colorlinks=true, citecolor=blue}  
\renewcommand{\Re}{\operatorname{Re}}

\DeclareGraphicsExtensions{.png,.pdf,.eps,.jpg}
\title[{\tiny On the impact of boundary conditions on weakly coupled thermoelastic wave model}]{
 \bf On the impact of boundary conditions on weakly coupled thermoelastic wave model}
\author{S. Nafiri} 
\thanks{The author would like to thank Professor E. Zuazua for fruitful discussions and remarks during his internship in DeustoTech Research Center of the University of Deusto in Bilbao.}

\address{Salem Nafiri \newline
Département de Mathématiques, informatique et géomatique. Ecole Hassania des Travaux Publics. Route d'El Jadida, B.P 8108, Oasis, Casablanca, Morocco.}
\email{nafirisalem@gmail.com/nafiri@ehtp.ac.ma}

\date{}
\begin{document}
\maketitle
\numberwithin{equation}{section}
\newtheorem{thm}{Theorem}[section]
\newtheorem{lem}{Lemma}[section]
\newtheorem{prop}{Proposition}[section]
\newtheorem{Def}{Definition}[section]
\newtheorem{rmk}{Remark}[section]
\newenvironment{dem}{\noindent\textbf{Proof.~}}{\hfill$\square$\bigbreak}

\begin{abstract}
The purpose of this paper is to demonstrate how different types of boundary conditions do not impact the asymptotic behaviour of the solutions of thermoelastic wave model. For an initial-boundary value problem associated with this system, we prove a global well-posedness result in a certain topology under appropriate regularity conditions on the data. Further, we show that under particular classes of boundary conditions, the energy associated to the system decays polynomially to zero and not exponentially.\\

{\bf Key words:}
thermoelastic structure, contraction semigroups, exponential decay, polynomial decay.
\end{abstract}
\section{\bf Introduction}
\label{sect1}
It is well known from experiment that the deformation of a body is inseparably connected with a change of its heat content and therefore with a change of the temperature distribution in the body. A deformation of a body which varies in time leads to temperature changes, and conversely. The internal energy of the body depends on both the temperature and the deformation. The science which deals with the investigation of the above coupled processes, is called thermoelasticity. The mathematical model for such (thermoelastic) bodies usually consists of a pair of coupled partial differential equations, one which models the elastic displacement (such as wave, beam or plate equation) and the other which models thermal diffusion. We refer to \cite[chap V, p.203]{Walker1980} and \cite{Nowacki1970} for references and discussion of modeling issues.

As far as linear thermoelastic systems are concerned, Dafermos \cite{Dafermos1968} was probably the first to investigate the asymptotic behavior of solutions to the following system 
\begin{equation}
\label{eq1.1}
 \left\{\begin{array}{lll}
                      u_{tt}(x,t)- \Delta u(x,t)+\gamma~\nabla\theta(x,t)=0\quad &in\;\Omega\times(0, \infty),\\
		      \theta_{t}(x,t)-\Delta\theta(x,t)-\gamma~div~u_{t}(x,t) =0\quad &in\;\Omega\times(0, \infty),\\
u(x,0)=u_{0}(x),\; u_{t}(x,0)=u_{1}(x),\; \theta(x,0)=\theta_{0}(x) \quad &on \;\Omega.
                     \end{array}
                   \right.
\end{equation}
In \cite{Dafermos1968}, he showed that the solutions to this thermoelastic system approach an equilibrium. However, the question of how fast do the solutions approach it remained open for some time. 
Hansen \cite{Hansen1992} succeeded in establishing an exponential decay estimate (for special boundary conditions) by using Fourier series expansion and a decoupling technique. Afterward, an essential progress was achieved in treating such (exponential) behavior for other boundary conditions, see the book by Liu \& Zheng \cite[Chap 2, p:23-25]{LZb1999} for a chronoligical treatment of the topic.

In this paper, we are concerned with the polynomial decay estimate of the energy associated to a (weakly) linear thermoelastic system, subject to various boundary conditions. We mean by (weakly) linear thermoelastic system, system \eqref{eq1.1} where the coupling terms ($\gamma\;\nabla\theta$ and $\gamma\; div\;u_{t}$) are replaced by ($\gamma\theta$ and $\gamma u_{t}$). With these coupling terms, the thermoelastic system considered is different to the one intensively treated in litterature, in the sense that its solutions do not decay exponentially to the equilibrium. A fact that has been shown by Khodja \& al. in \cite{KBT1997}. Our purpose is to study the effect of boundary conditions on the rate of decay of solutions associated to the system in consideration. 

Let $u(x,t)$ be the displacement (longitudinal or transverse, depending upon the application) at position $x$ along a bounded smooth domain $\Omega\subset\mathbb{R}^{n}$ and time $t$ , and $\theta(x,t)$ be the temperature deviation from the reference temperature at position $x$ and time $t$. The coupling constant $\gamma$ is generally small in comparison to 1 and is a measure of the mechanical-thermal coupling present in the system (see \cite{Hansen1992}, \cite{JLions1988}).  Then the initial value problem governed by $u$ and $\theta$ is described as follows
\begin{equation}
\label{eq1.2}
 \left\{\begin{array}{lll}
                      u_{tt}(x,t)- \Delta u(x,t)+\gamma \theta(x,t)=0\quad &in\;\Omega\times(0, \infty),\\
		      \theta_{t}(x,t)-\Delta\theta(x,t)-\gamma u_{t}(x,t) =0\quad &in\;\Omega\times(0, \infty),\\
u(x,0)=u_{0}(x),\; u_{t}(x,0)=u_{1}(x),\; \theta(x,0)=\theta_{0}(x) \quad &on \;\Omega,
                     \end{array}
                   \right.
\end{equation}

The set of boundary conditions (\textbf{B.C}) that we considered here are the following

\begin{equation}
\label{DD}
 u=0=\theta \qquad \text{on}\; \partial\Omega\times (0, \infty),\qquad\text{\textbf{(Dirichlet-Dirichlet B.C)}} \tag{DD}
\end{equation}
\begin{equation}
\label{DN}
 u=0=\frac{\partial\theta}{\partial n} \;\;\quad \text{on}\; \partial\Omega\times (0, \infty),\qquad\text{\textbf{(Dirichlet-Neumann B.C)}} \tag{DN}
\end{equation}
\begin{equation}
\label{ND}
\frac{\partial u}{\partial n}=0=\theta \;\;\quad \text{on}\; \partial\Omega\times (0, \infty),\qquad\text{\textbf{(Neumann-Dirichlet B.C)}}\tag{ND}
\end{equation}

\begin{equation}
\label{NN}
\frac{\partial u}{\partial n}=0=\frac{\partial\theta}{\partial n} \quad \text{on}\; \partial\Omega\times (0, \infty),\qquad\text{\textbf{(Neumann-Neumann B.C)}}\tag{NN}
\end{equation}
where $n$ is the unit outward normal to $\Omega$.

The energy associated to system \eqref{eq1.2} is given by
\begin{equation}
\label{eq1.3}
E(t)=\int_{\Omega}|\nabla u(x,t)|^{2}+|u_{t}(x,t)|^{2}+|\theta(x,t)|^{2}dx.
\end{equation}

The outline of this paper is the following. In section 2, we show a well posedness result of system \eqref{eq1.2} in appropriate energy space. In section 3, we show that the semigroup associated to system \eqref{eq1.2} is strongly stable. Finally, in section 4, we prove that the semigroup associated to system \eqref{eq1.2} does not decays exponentially to zero, but polynomially.
\section{Well-posedness of thermoelastic system in abstract setting}
Let $H$ be a complex Hilbert space with the inner product $\langle\cdot,\cdot\rangle$ and the induced norm $\|\cdot\|$. We consider in the following the thermoelastic system \eqref{eq1.2} in abstract setting:
\begin{equation*}
\left\{
\begin{array}{lll}
u_{tt}= -A_Ou -\gamma\theta\\
\theta_{t}=-A_{O'}\theta +\gamma u_{t}\\
u(0)=u_0,\quad u_t(0)=v_0,\quad \theta(0)=\theta_0
\end{array}
\right.
\end{equation*}
where $A_O$ and $A_{O'}$, $O,O'\in\{D,N\}$, are self-adjoint, strictly positive definite (unbounded) operators with compact resolvent on a complex Hilbert space $H$. We define the state (energy) space 
$$
\mathcal{H}_{O}=D(A^{\frac{1}{2}}_{O})\times H\times H,\qquad O\in\{D,N\}.
$$

Any element in $\mathcal{H}_{O}$ is denoted by $U=(u,v,\theta)^T$. Introduce
\begin{eqnarray}
\label{sp}
\langle U_1,U_2\rangle_{\mathcal{H}_{O}}=\langle A^{\frac{1}{2}}_Ou_1,A^{\frac{1}{2}}_O u_2\rangle+\langle v_1,v_2\rangle+\langle\theta_1,\theta_2\rangle
\end{eqnarray}
for all $U_i=(u,v,\theta)^T\in\mathcal{H}_{O},\; i=1,2,\; O\in\{D,N\}$ and the induced norm
\begin{eqnarray*}
\|(u,v,\theta)\|^2_{\mathcal{H}_{O}}=\|A^{\frac{1}{2}}_O u\|^2+\|v\|^2+\|\theta\|^2.
\end{eqnarray*}

By denoting the variable (velocity) $v=u_{t}$ and $U_{0}=(u_{0},u_{1},\theta_{0})^{T}$, system \eqref{eq1.2}-(\textbf{B.C}) can be written in the following abstract first-order evolution equation on the space $\mathcal{H}_{O}$,
\begin{equation}
\label{eq2.1}
 \left\{\begin{array}{ll}
 \frac{dU(t)}{dt}=\mathcal{A}_{OO'}U(t),\quad t\geq 0\\
 U(0)=U_{0},\\
\end{array}
\right.
\end{equation}
where the operator $\mathcal{A}_{OO'}:D(\mathcal{A}_{OO'})\subseteq\mathcal{H}_{O}\rightarrow \mathcal{H}_{O}$ is defined by
\begin{equation}
\label{eq2.3}
\mathcal{A}_{OO'}=\left(\begin{matrix}
                    0 & 1 & 0\\
	       -A_{O} & 0 & -\gamma 1\\
		    0 & \gamma 1 & -A_{O'}
                   \end{matrix}
\right),\quad O,O'\in\{D,N\}.
\end{equation}

\[
 \mathcal{A}_{OO'}\in\Big\{\mathcal{A}_{DD},\mathcal{A}_{DN},\mathcal{A}_{ND},\mathcal{A}_{NN}\Big\}
\]
with the domain
$$
D(\mathcal{A}_{OO'})=D(A_{O})\times D(A^{\frac{1}{2}}_{O})\times D(A_{O'}),\quad O,O'\in\{D,N\}
$$
and
\[
\mathcal{A}_{OO'}(u,v,\theta)=(v,-A_{O} u-\gamma\theta,-A_{O'}\theta+\gamma v),\quad O,O'\in\{D,N\}.
\]
Here we have used the notation $A_{O}$ and $A_{O'}$ to distinguish between the Dirichlet Laplacian and the Neumann Laplacian. The operators $A_{D}$ and $A_{N}$ of $L^{2}(\Omega)$ are defined by 
\begin{align*}
A_{D}\phi_m=-\Delta \phi_m=\lambda_m \phi_m,&\qquad D(A_{D})=H^{2}(\Omega)\cap H_{0}^{1}(\Omega)\\
A_{N}\xi_m:=-\Delta \xi_m=\mu_m\xi_m,&\qquad D(A_{N})=\{w\in H^{2}(\Omega)/ \frac{\partial w}{\partial n}\mid_{\partial\Omega}=0\}
\end{align*}
and by \cite{G1967}, we have
\[
D(A_D^\frac{1}{2})=H_0^1(\Omega),\;\text{ and }D(A_N^\frac{1}{2})=H^1(\Omega).
\]
$(\lambda_m,\phi_m)$ and $(\mu_m,\xi_m)$ are the eigenvalues and eigenvectors of the operators $A_D$ and $A_N$ respectively, with
$$
\lim\limits_{m \rightarrow +\infty}\lambda_m=\lim\limits_{t \rightarrow +\infty}\mu_m=+\infty,\quad \|\phi_m\|_H=\|\xi_m\|_H=1.
$$
Th following theorem states that system \eqref{eq1.2} subject to (\textbf{B.C}) is  well-posed in the Hilbert space $\mathcal{H}_O$, $O\in\{D,N\}$. To show the well-posedness, we will rely on the following lemma.

\begin{thm}
\label{thm2.1}
For all $O,O'\in\{D,N\}$ the family of operators $\mathcal{A}_{OO'}$, generates a contraction semigroup $T_{OO'}(\cdot)$ on the Hilbert space $\mathcal{H}_O$.
\end{thm}

\begin{dem}
To do this, we need to show that $\mathcal{A}_{OO'}$, $O,O'\in\{D,N\}$, is m-dissipative  on the Hilbert space $\mathcal{H}_O$, $O\in\{D,N\}$. 
\begin{itemize}
\item[(a)] Indeed, for any $U\in D(\mathcal{A}_{OO'})$, $O,O'\in\{D,N\}$, by using the definition \eqref{sp}, we have
\begin{align}
\label{eq2.4}
\langle\mathcal{A}_{OO'}U,U\rangle_{\mathcal{H}_O} &=
\langle (v,-A_Ou-\gamma\theta,\gamma v-A_O'\theta), (u,v,\theta)\rangle_{\mathcal{H}_O}\nonumber\\
&=\langle A_{O}^{\frac{1}{2}}v,A_{O}^{\frac{1}{2}}u\rangle-\langle A_O u ,v \rangle -\gamma\langle\theta,v \rangle+\gamma\langle v,\theta\rangle -\langle A_O' \theta,\theta\rangle\nonumber\\
\Re \langle\mathcal{A}_{OO'}U,U\rangle_{\mathcal{H}_O} &=-\|A_{O'}^{\frac{1}{2}}\theta\|^2 \leqslant 0.
\end{align}
This implies that $\mathcal{A}_{OO'}$ is a dissipative operator. Here Re is used to denote the real part of a complex number.\\
\item[(b)] Let's show now that $[\mu I-\mathcal{A}_{OO'}]$, $O,O'\in\{D,N\}$, is surjective. That is, for all $F=(f,g,h)^{T}\in\mathcal{H}_O$, $O\in\{D,N\}$, there exist $U=(u,v,\theta)^T\in D(\mathcal{A}_{OO'})$, $O,O'\in\{D,N\}$, satisfying
\begin{align}
\label{eq2.5}
[I-\mathcal{A}_{OO'}]U=F,
\end{align}
i.e.,
\begin{align}
u-v &=f,\label{eq2.6}\\
v+A_{O}u+\gamma\theta &=g,\label{eq2.7}\\
\theta+A_{O'}\theta-\gamma v &=h.\label{eq2.8}
\end{align}
Equation (2.6) implies that $v=u-f$ in $D(A^\frac{1}{2})$. Replacing this in (2.7) and (2.8), we obtain
\begin{align*}
(I+A_O)u+\gamma\theta &=f+g\\
(I+A_O')\theta-\gamma u &=h-\gamma f
\end{align*}
which is equivalent in matrix form to
\begin{align*}
\left(\begin{matrix}
          I+A_O  & \gamma \\
	       -\gamma & I+A_O'  
                   \end{matrix}\right)
\left(\begin{matrix}
          u \\
	    \theta  
\end{matrix}\right)=\left(\begin{matrix}
          f+g \\
	    h-\gamma f  
\end{matrix}\right).
\end{align*}
We set
\begin{align*}
L\equiv\left(\begin{matrix}
          I+A_O  & \gamma \\
	       -\gamma & I+A_O'  
                   \end{matrix}\right)
=\left(\begin{matrix}
          I+A_O  & 0 \\
	       0 & I+A_O'  
                   \end{matrix}\right)
+\left(\begin{matrix}
          0  & \gamma \\
	       -\gamma & 0  
                   \end{matrix}\right)\equiv M+N.
\end{align*}
Clearly $M$ is a generator of contraction semigroup and since $N$ is dissipative, $L$ is also a generator of contraction semigroup applying the Bounded Perturbation Theorem and/or perturbation result obtained by Desch and Schappacher in \cite[p. 335]{DS1984}, \cite[Thm 2.7, p.173]{EN2000}, cf. also \cite[Remark 7, p. 256]{KBT1997}.
\end{itemize}
The conclusion immediately follows from the Lumer-Phillips theorem \cite[Corollary II.3.20]{EN2000}.
\end{dem}
\begin{rmk}
We could also use the fact that for all $O,O'\in\{D,N\}$, $\mathcal{A}_{OO'}$ and $\mathcal{A}^*_{OO'}$ are dissipative and conclude by the Lumer-Phillips theorem.
\end{rmk}
\begin{rmk}
If we replace system \eqref{eq1.2}, by the following one
\begin{equation}
\label{eq2.59}
\left\{
\begin{array}{lll}
u_{tt}= -A_Ou -\alpha\theta,\quad O\in\{D,N\}\\
\theta_{t}=-A_{O'}\theta -\beta u_{t},\quad O'\in\{D,N\}\\
u(0)=u_0,\quad u_t(0)=v_0,\quad \theta(0)=\theta_0
\end{array}
\right.
\end{equation}
where $\alpha$ and $\beta$ are arbitrary constants in $\mathbb{R}$, we still can show in a same way well posedness of \eqref{eq2.59}-(\textbf{B.C}), provided that we deal with a real Hilbert space $H$. Indeed, for all $O,O'\in\{D,N\}$, we have\\
\begin{align}
\label{eq2.60}
\langle\mathcal{A}_{OO'}U,U\rangle_{\mathcal{H}_O} &=\langle A_O^{\frac{1}{2}}v, A_O^{\frac{1}{2}}u\rangle-\langle A_O u,v\rangle-\alpha\langle\theta,v\rangle+\langle A_{O'}\theta,\theta\rangle-\beta\langle v,\theta\rangle\nonumber\\
					&=-(\alpha+\beta)\langle\theta,v\rangle-\|A_{O'}^{\frac{1}{2}}\theta\|^{2} \nonumber\\
					&\leqslant |\alpha+\beta|\langle\theta, v\rangle-\|A_{O'}^{\frac{1}{2}}\theta\|^{2} \nonumber\\
					&\leqslant \frac{(\alpha+\beta)^{2}}{2}\|v\|^{2}+\frac{1}{2}\|\theta\|^{2}-\|A_{O'}^{\frac{1}{2}}\theta\|^{2} \nonumber\\
					&\leqslant\frac{\max\{(\alpha+\beta)^{2},1\}}{2}\|U\|^{2}_{\mathcal{H}_O}-\|A_{O'}^{\frac{1}{2}}\theta\|^{2}\nonumber\\
					&=\frac{C}{2}\|U\|^{2}_{\mathcal{H}_O}-\|A_{O'}^{\frac{1}{2}}\theta\|^{2}.
\end{align}
Introducing $\mathcal{B}_{OO'}=-\frac{C}{2}I+\mathcal{A}_{OO'}$ and following the proof of Theorem  \ref{thm2.1}, we can prove that $\mathcal{A}_{OO'}$, $O,O'\in\{D,N\}$ generates a contraction semigroup on $\mathcal{H}_O$, $O\in\{D,N\}$. However, to show the asymptotical behaviour of \eqref{eq2.59}, we need $\alpha$ and $\beta$ to be nonzero real numbers with the same sign, see \cite{R1992,T2010}.
\end{rmk}

\section{Asymptotic behaviour}
\subsection{Strong stability}
To show that the semigroup associated to system \eqref{eq1.2} is strongly stable, we introduce here the notions of stability that we encounter in this work and we consider the following energy spaces
\begin{align*}
X_{DD}&=D(A_D^\frac{1}{2})\times H\times H\\
X_{DN}&=D(A_D^\frac{1}{2})\times H\times H_N\\
X_{ND}&=D(A_N^\frac{1}{2})\times H_N\times H\\
X_{NN}&=D(A_N^\frac{1}{2})\times H_N\times H_N
\end{align*}
where
\[
H_{N}=\{f\in H: \langle f,1\rangle_H=0\}.
\]
\begin{Def}
Assume that $\mathcal{A}$ is the generator of a $C_0-$semigroup of contractions $(T(t))_{t\geqslant 0}$ on a Hilbert space $\mathcal{H}$. The $C_0-$semigroup $T(\cdot)$ is said to be
\begin{enumerate}
\item[1.] strongly stable if
\begin{align*}
\lim\limits_{t \rightarrow +\infty} \|T(t)U_{0}\|_{X}=0,\quad U_{0}\in\mathcal{H};
\end{align*}
\item[2.] exponentially (or uniformly) stable with decay rate $\omega > 0$ if there exists a
constant $M\geqslant 1$ such that 
$$
\|T(t)U_0\|_\mathcal{H}\leqslant Me^{-\omega t}\|U_0\|_\mathcal{H},\quad t>0,\; U_0\in\mathcal{H};
$$
\item[3.] polynomially stable (of order $\alpha>0$) if it is bounded, if $i\mathbb{R}\subset\rho(\mathcal{A})$ and if
$$
\|T(t)U_0\|_\mathcal{H}\leqslant Ct^{-\alpha}\|\mathcal{A}U_0\|_\mathcal{H},\quad t>0,\; U_0\in D(\mathcal{A})
$$
for some positive constant $C$.\\
In that case, one says that solutions of \eqref{eq2.1} decay at a rate $t^{-\alpha}$. The $C_0-$semigroup $T(\cdot)$ is said to be polynomially stable with optimal decay rate $t^{-\alpha}$ (with $\alpha>0$) if it is polynomially stable with decay rate $t^{-\alpha}$ and, for any $\epsilon>0$ small enough, there exists solutions of \eqref{eq2.1} which do not decay at a rate $t^{-(\alpha-\epsilon)}$.
\end{enumerate}

\end{Def}

To show the strong stability of the $C_0-$semigroup $T_{OO'}(\cdot)$, $O,O'\in\{D,N\}$. We will rely on the following result obtained first by Benchimol in \cite{B1978} and then by Arendt and Batty in \cite{AB1988}.
\begin{lem}[\textbf{Benchimol, Arendt, Batty in} \cite{AB1988,B1978}]
Assume that $\mathcal{A}_{OO'}$ is the generator of a $C_0-$semigroup of contractions $T_{OO'}(t)$ on a Hilbert space $\mathcal{H}_{O}$, $O,O'\in\{D,N\}$. If $i\mathbb{R}\subset\rho
(\mathcal{A}_{OO'})$, $O,O'\in\{D,N\}$, then the $C_0-$semigroup $T_{OO'}(t)$ is strongly stable on $\mathcal{H}_{O}$, $O,O'\in\{D,N\}$.
\end{lem}

\begin{lem}
For all $O,O'\in\{D,N\}$, the family of operators $\mathcal{A}_{OO'}:D(\mathcal{A}_{OO'})\subseteq X_{OO'}\rightarrow X_{OO'}$, satisfies $i\mathbb{R}\subset\rho(\mathcal{A}_{OO'})$, .
\end{lem}
\begin{dem}
By contradiction argument, let $0\ne U =(u,v,\theta)^T\in D(\mathcal{A}_{OO'})$, $\beta\in\mathbb{R}$ such that
\begin{equation}
\label{eq3.1}
\mathcal{A}_{OO'}U=i\beta U.
\end{equation}
Our goal is to find a contradiction by proving that $U = 0$. 

If $\beta=0$, we have $\mathcal{A}_{OO'}U=0$, that is
\begin{align*}
v=0\\
-A_O u-\gamma\theta=0\\
\gamma v -A_{O'}\theta=0.
\end{align*}
\textbf{Case 1: Dirichlet-Dirichlet B.C.}
This case is obvious since $0\in\rho(A_D)$.\\
\textbf{Case 2: Dirichlet-Neumann B.C.}
$v=0$ implies $\theta=0$, since $\theta\in H_N$ and then $u=0$.\\
\textbf{Case 3: Neumann-Dirichlet B.C.}
$v=0$ implies $\theta=0$ and then $u=0$, because $u\in H_N$.\\
\textbf{Case 4: Neumann-Neumann B.C.}
$v=0$ implies $\theta=u=0$, since $\theta,u\in H_N$.\\ 
In all cases, we deduce that $U=0$. In the following we assume that $\beta\ne 0$. Taking the real part of the inner product in $X_{OO'}$ of $\mathcal{A}_{OO'}U$ and $U$, we obtain
\begin{align*}
0=\Re(i\beta\|U\|^{2}_{\mathcal{H}_O})=\Re\langle\mathcal{A}_{OO'}U,U\rangle_{\mathcal{H}_O}=-\|A^{\frac{1}{2}}_{O'}\theta\|^{2}.
\end{align*}
\textbf{Case 1: Dirichlet-Dirichlet B.C.}
Dissipativity implies that $\theta=0$, replacing in \eqref{eq3.1}, we obtain $v=0$ and then $u=0$.\\
\textbf{Case 2: Dirichlet-Neumann B.C.}
Dissipativity with $\theta\in H_N$, replacing in \eqref{eq3.1} implies $\theta=0$ and then $v=u=0$.\\
\textbf{Case 3: Neumann-Dirichlet B.C.}
Dissipativity implies that $\theta=0$, replacing in \eqref{eq3.1}, we obtain $v=u=0$.\\
\textbf{Case 4: Neumann-Neumann B.C.}
Dissipativity with $\theta\in H_N$, replacing in \eqref{eq3.1} implies $\theta=0$ we obtain $v=u=0$.

We conclude that $U=0$ and the desired contradiction is proved.
\end{dem}

\begin{thm}
The semigroup $T_{OO'}(t)$ is strongly stable in the energy space $\mathcal{H}_O$, $O,O'\in\{D,N\}$. In
other words
\begin{align*}
\lim\limits_{t \rightarrow +\infty} \|T_{OO'}(t)U_{0}\|_{X}=0,\quad U_{0}\in \mathcal{H}_O,\;O,O'\in\{D,N\}.
\end{align*}
\end{thm}

\begin{dem}
The proof follows from Lemma 3.1 and Lemma 3.2.
\end{dem}
\subsection{Lack of exponential decay}
Our starting point in this section is to show the lack of exponential decay is of the solutions of the system \eqref{eq2.1}. That is, the resolvent operator is not uniformly bounded. To do this, we use the following theorem (frequency domain method) which has been proved independently by Gearhart \cite{G1978}, Pr\"uss \cite{P1984} and Huang \cite{Hu1985}.
\begin{thm}[\textbf{Gearhart, Pr\"uss, Huang in} \cite{G1978,P1984,Hu1985}]
Let $T(t)$ be a $C_{0}-$semigroup of contractions of linear operators on
Hilbert space with infinitesimal generator $\mathcal{A}$. Then $T(t)$ is exponentially stable if and
only if
\begin{itemize}
\item[(a)] $i\mathbb{R}\subset \rho(\mathcal{A})$\\
\item[(b)] $\underset{|\beta|\to\infty}{sup}\|(i\beta I-\mathcal{A})^{-1}\|_{\mathcal{L}(\mathcal{H})}=O(1)$.
\end{itemize}
\end{thm}
We are now in conditions to show the main result of this section.
\begin{thm}
For all $O,O'\in\{D,N\}$, the semigroup $T_{OO'}(\cdot)$ is not exponentially stable.
\end{thm}

\begin{dem}
To prove Theorem 4.2, it is sufficient to show that the condition (b) of the Theorem 4.1 does not hold. To do this, it is sufficient to show the existence of sequences $F_n\in\mathcal{H}_O$, $O\in\{D,N\}$ and $\beta_n\in\mathbb{R}$ such that $(F_n)_{n\in\mathbb{N}}$ is bounded, $|\beta|\rightarrow\infty$ and $\|(i\beta I-\mathcal{A}_{OO'})^{-1}F_n\|\rightarrow\infty$ when $n\to\infty$.

Let us take $\lambda\in\mathbb{R}$. Then for any $F_n=(f_n,g_n,h_n)^T\in\mathcal{H}$ there exists only one $U_n=(i\lambda I-\mathcal{A}_{OO'})^{-1}F_n=(u_n,v_n,\theta_n)^T\in D(\mathcal{A}_{OO'})$ solution of the resolvent equation
\begin{align*}
i\lambda U_n-\mathcal{A}_{OO'}U_n=F_n.
\end{align*}
\textbf{Case 1: Dirichlet-Dirichlet B.C.}
To simplify the notation we will omit the sub index $m$. The equation reads
\begin{align*}
i\lambda u-v=f\\
A_D u+i\lambda v+\gamma\theta=g\\
-\gamma v+(i\lambda+A_D)\theta=h.
\end{align*}
Taking $f=0=h$ and $g=\phi_m$, the solutions of the previous system must be of the form $u=a\phi_m$, $v=i\lambda u=i\lambda a\phi_m$, and $\theta=b\phi_m$, where $a$ and $b$ verify
\begin{align*}
a(\lambda_m-\lambda^2)+b\gamma=1\\
i\lambda a\gamma+b(i\lambda+\lambda_m)=0.
\end{align*}
Now, choosing $\lambda=\sqrt{\lambda_m}$, we obtain
$b=\frac{1}{\gamma}$ and $a=\frac{1}{\gamma^2}(1-i\sqrt{\lambda_m})$.\\
Recalling that $u=a\phi_m=\frac{1}{\gamma^2}(1-i\sqrt{\lambda_m})\phi_m$, we get
$$
\|u\|^2=\frac{1}{\gamma^4}(\lambda+\lambda^2_m).
$$
Therefore we have
$\lim\limits_{m \rightarrow +\infty}\|U_m\|^2_{\mathcal{H}_D}\geqslant\lim\limits_{m \rightarrow +\infty}\|u\|^2=\infty$.\\
\textbf{Case 2: Dirichlet-Neumann B.C.}
We proceed in the same way, we take $f=0=h$, $g=\phi_m$ and we look for solutions in the form $u=a\phi_m$, $v=i\lambda u=i\lambda a\phi_m$, and $\theta=b\xi_m$. Simple calculations give that $b=\frac{1}{\gamma}$ and $a=\frac{1}{\gamma^2}\sqrt{1+\frac{\lambda_m^2}{\mu_m}}$. Since $\|u\|^2=\frac{1}{\gamma^4}(\mu_m+\lambda_m^2)$, we obtain
$$
\lim\limits_{m \rightarrow +\infty}\|U_m\|^2_{\mathcal{H}_D}\geqslant\lim\limits_{m \rightarrow +\infty}\|u\|^2=\infty.
$$
\textbf{Case 3: Neumann-Dirichlet B.C.}
We proceed in the same way, we take $f=0=h$, $g=\xi_m$ and we look for solutions in the form $u=a\xi_m$, $v=i\lambda u=i\lambda a\xi_m$, and $\theta=b\phi_m$. Simple calculations give that $b=\frac{1}{\gamma}$ and $a=\frac{1}{\gamma^2}\sqrt{1+\frac{\mu_m^2}{\lambda_m}}$. Since $\|u\|^2=\frac{1}{\gamma^4}(\lambda_m+\mu_m^2)$, we obtain
$$
\lim\limits_{m \rightarrow +\infty}\|U_m\|^2_{\mathcal{H}_D}\geqslant\lim\limits_{m \rightarrow +\infty}\|u\|^2=\infty.
$$
\textbf{Case 4: Neumann-Neumann B.C.}
We proceed in the same way, we take $f=0=h$, $g=\xi_m$ and we look for solutions in the form $u=a\xi_m$, $v=i\lambda u=i\lambda a\xi_m$, and $\theta=b\xi_m$. Simple calculations give that $b=\frac{1}{\gamma}$ and $a=\frac{1}{\gamma^2}(1-i\sqrt{\mu_m})$. Since $\|u\|^2=\frac{1}{\gamma^4}(\mu_m+\mu_m^2)$, we obtain
$$
\lim\limits_{m \rightarrow +\infty}\|U_m\|^2_{\mathcal{H}_D}\geqslant\lim\limits_{m \rightarrow +\infty}\|u\|^2=\infty
$$
which completes the proof.
\end{dem}
\subsection{Polynomial decay}
In the case where $T_{OO'}(\cdot),\;O,O'\in\{D,N\}$ is not exponentially stable, we look for a polynomial decay rate. In this section, we use the frequency domain approach method to show the polynomial stability of $T_{OO'}(\cdot),\;O,O'\in\{D,N\}$ associated with the weakly coupled thermoelastic wave model \eqref{eq1.2} subject to (\textbf{B.C}). The frequency domain approach method has been obtained by Batty in \cite{BD2008}, Borichev and Tomilov in \cite{BT2010}, Liu and Rao in \cite{LR2005}, B\' atkai and al. in \cite{BEPS2006}.
\begin{thm}[\textbf{Borichev-Tomilov in \cite{BT2010}}]
Assume that $\mathcal{A}$ is the generator of a strongly continuous $C_0-$semigroup of contractions $T(\cdot)$ on a Hilbert space $\mathcal{H}$. If $i\mathbb{R}\subset\rho (\mathcal{A})$, then for a fixed $\alpha>0$ the following conditions are equivalent
\begin{itemize}
\item[(i)] $\underset{\lambda\in\mathbb{R}}{sup}\|(i\lambda I-\mathcal{A})^{-1}\|_{\mathcal{L}(\mathcal{H})}=O(|\lambda|^{\alpha})$
\item[(ii)] $\|T(t)U_0\|\leqslant C t^{-\frac{1}{\alpha}}\|U_0|_{D(\mathcal{A})},\;t>0,\;U_0\in D(\mathcal{A}),\;\text{for some } C>0$.
\end{itemize}
\end{thm}

\begin{thm}
For all $O,O'\in\{D,N\}$, the semigroup $T_{OO'}(t)$, is polynomially stable (of order $\alpha=2$) on the Hilbert space $\mathcal{H}_O$.
\end{thm}

\begin{dem}
For any $\lambda\in\mathbb{R}$, and any $U\equiv (u,v,\theta)^T\in D(\mathcal{A}_{OO'})$,

\begin{align*}
(i\lambda I-\mathcal{A})U=\left(\begin{matrix}
                    i\lambda & -1 & 0\\
	       A_O & i\lambda & \gamma 1\\
		    0 & -\gamma 1 & i\lambda+A_{O'}
                   \end{matrix}
\right)\left(\begin{matrix}
           u\\
	       v\\
		   \theta
                   \end{matrix}
\right)=\left(\begin{matrix}
                    i\lambda u-v\\
	       A_O u+i\lambda v+\gamma\theta\\
		    -\gamma v+(i\lambda+A_{O'})\theta
                   \end{matrix}
\right).
\end{align*}
Our proof  will be based on a contradiction argument. Suppose that the following is not true.
\[
\underset{\lambda\in\mathbb{R}}{sup}|\lambda|^{-2}\|(i\lambda I-\mathcal{A}_{OO'})^{-1}\|_{\mathcal{L}(\mathcal{H})}<\infty.
\]
Then there exists a sequence $\{(\lambda_n,U_n)|n\geqslant 1\}\subseteq\mathbb{R}\times D(\mathcal{A})$, with $U_n\equiv(u_n,v_n,\theta_n)^T$, and
\begin{align*}
\left\{\begin{array}{ll}
 \lim\limits_{n \to +\infty}|\lambda_n|=\infty\\
 \|U_n\|_\mathcal{H}=\|A_O^\frac{1}{2}u_n\|^2+\|v_n\|^2+\|\theta_n\|^2=1,\;n\geqslant 1,
\end{array}
\right.
\end{align*}
such that
\[
\lim\limits_{n \to +\infty}|\lambda|^{2}\|(i\lambda_n-\mathcal{A}_{OO'})U_n\|_\mathcal{H}=0,
\]
i.e.
\begin{equation}
\label{eq4.1}
\lambda_n^{2}(i\lambda_n A_O^\frac{1}{2}u_n-A_O^\frac{1}{2}v_n)=o(1)
\end{equation}
\begin{equation}
\label{eq4.2}
\lambda_n^{2}(i\lambda_n v_n+A_O u_n+\gamma\theta_n)=o(1)
\end{equation}
\begin{equation}
\label{eq4.3}
\lambda_n^{2}(i\lambda_n\theta_n-\gamma v_n+A_{O'}\theta_n)=o(1).
\end{equation}
Hereafter $o(1)$ stands for a vector in $H$ (or a quantity in $\mathbb{R}$) which goes to zero as $n\to\infty$. The advantage of using such a notation is that the previous system can be regarded as a system of equations, which will be convenient below.

Our goal is to obtain $\|U_n\|=o(1)$, thus a contradiction. By the dissipativeness of the operator $\mathcal{A}$
\begin{equation}
\label{eq3.5}
-|\lambda_n|^2\|A_{O'}^\frac{1}{2}\theta_n\|^2=Re\langle |\lambda_n|^2(i\lambda_n I-\mathcal{A}_{OO'})U_n,U_n\rangle_\mathcal{H}=o(1).
\end{equation}
\textbf{Case 1: Dirichlet-Dirichlet B.C.}
From \eqref{eq3.5}, we have $\|A_{D}^\frac{1}{2}\theta_n\|=o(1)$ and since $0\in\rho(A_D)$ it follows that $\theta=o(1)$.\\
Taking the inner product of $\eqref{eq4.3}$ with $v_n$ in $H$ yields
\[
i\langle \lambda_n\theta_n,v_n\rangle-\gamma\|v_n\|^2+\langle \lambda_n A_D^\frac{1}{2}\theta_n,\frac{1}{\lambda_n} A_D^\frac{1}{2}v_n\rangle=o(1).
\]
It follows from the dissipativeness property and $\eqref{eq4.1}$ that $\frac{1}{\lambda_n} A^\frac{1}{2}v_n$ is bounded in $H$. Thus the third term in the previous expression converges to zero. The first term in the previous expression also converges to zero due to $\theta=0$ and the boundedness of $v_n$. We have obtained $v_n=o(1)$.\\
On the other hand, using the sum of the inner product of $\eqref{eq4.1}$ with $u_n$ in $D(A^\frac{1}{2})$ and the inner product of $\eqref{eq4.2}$ with $-v_n$ in $H$,
\[
\|A_D^\frac{1}{2} u_n\|^2-\|v_n\|^2=o(1).
\]
\textbf{Case 2: Dirichlet-Neumann B.C.}
From \eqref{eq3.5}, we have $\|A_{N}^\frac{1}{2}\theta_n\|=o(1)$ and since $\theta\in H_N$, then $\theta_n=o(1)$.\\
Taking the inner product of $\eqref{eq4.3}$ with $v_n$ in $H$ yields
\[
i\langle \lambda_n\theta_n,v_n\rangle-\gamma\|v_n\|^2+\langle \lambda_n A_N^\frac{1}{2}\theta_n,\frac{1}{\lambda_n}A_N^\frac{1}{2}v_n\rangle=o(1).
\]
It follows from the dissipativeness property and $\eqref{eq4.1}$ that $\frac{1}{\lambda_n}A_N^\frac{1}{2}v_n$ is bounded in $H$. Thus the third term in the previous expression converges to zero. The first term in the previous expression also converges to zero due to $\theta=0$ and the boundedness of $v_n$. We have obtained $v_n=o(1)$.\\
On the other hand, using the sum of the inner product of $\eqref{eq4.1}$ with $u_n$ in $D(A_D^\frac{1}{2})$ and the inner product of $\eqref{eq4.2}$ with $-v_n$ in $H$,
\[
\|A_D^\frac{1}{2} u_n\|^2-\|v_n\|^2=o(1).
\]
\textbf{Case 3: Neumann-Dirichlet B.C.}
From \eqref{eq3.5}, we have $\|A_{D}^\frac{1}{2}\theta_n\|=o(1)$ and since $0\in\rho(A_D)$ it follows that $\theta=o(1)$.\\
Taking the inner product of $\eqref{eq4.3}$ with $v_n$ in $H$ yields
\[
i\langle \lambda_n\theta_n,v_n\rangle-\gamma\|v_n\|^2+\langle \lambda_n A_D^\frac{1}{2}\theta_n,\frac{1}{\lambda_n}A_D^\frac{1}{2}v_n\rangle=o(1).
\]
It follows from the dissipativeness property and $\eqref{eq4.1}$ that $\frac{1}{\lambda_n}A_N^\frac{1}{2}v_n$ is bounded in $H$. Thus the third term in the previous expression converges to zero. The first term in the previous expression also converges to zero due to $\theta=0$ and the boundedness of $v_n$. We have obtained $v_n=o(1)$.\\
On the other hand, using the sum of the inner product of $\eqref{eq4.1}$ with $u_n$ in $D(A_N^\frac{1}{2})$ and the inner product of $\eqref{eq4.2}$ with $-v_n$ in $H$,
\[
\|A_N^\frac{1}{2} u_n\|^2-\|v_n\|^2=o(1).
\]
\textbf{Case 4: Neumann-Neumann B.C.}
From \eqref{eq3.5}, we have $\|A_{N}^\frac{1}{2}\theta_n\|=o(1)$ and since $\theta\in H_N$ it follows that $\theta=o(1)$.\\
Taking the inner product of $\eqref{eq4.3}$ with $v_n$ in $H$ yields
\[
i\langle \lambda_n\theta_n,v_n\rangle-\gamma\|v_n\|^2+\langle \lambda_n A_N^\frac{1}{2}\theta_n,\frac{1}{\lambda_n}A_N^\frac{1}{2}v_n\rangle=o(1).
\]
It follows from the dissipativeness property and $\eqref{eq4.1}$ that $\frac{1}{\lambda_n}A_N^\frac{1}{2}v_n$ is bounded in $H$. Thus the third term in the previous expression converges to zero. The first term in the previous expression also converges to zero due to $\theta=0$ and the boundedness of $v_n$. We have obtained $v_n=o(1)$.\\
On the other hand, using the sum of the inner product of $\eqref{eq4.1}$ with $u_n$ in $D(A_N^\frac{1}{2})$ and the inner product of $\eqref{eq4.2}$ with $-v_n$ in $H$,
\[
\|A_N^\frac{1}{2} u_n\|^2-\|v_n\|^2=o(1).
\]

We have the promised contradiction.

\end{dem}


\end{document}